\documentclass{article}

\usepackage[nonatbib, final]{nips_2016}

\usepackage[utf8]{inputenc} 
\usepackage[T1]{fontenc}    
\usepackage{hyperref}       
\usepackage{url}            
\usepackage{booktabs}       
\usepackage{amsfonts}       
\usepackage{microtype}      
\usepackage{graphicx}
\usepackage{algorithm}
\usepackage{algorithmic}

\newtheorem{theorem}{Theorem}
\newtheorem{lemma}[theorem]{Lemma}

\newtheorem{assumption}[theorem]{Assumption}

\title{Infinite dimensional optimistic optimisation with applications on physical systems}

\author{
    Muhammad F. Kasim, Peter A. Norreys \\
  Department of Physics, University of Oxford\\
  Oxford, OX1 3RH \\
  \texttt{\{m.kasim1, peter.norreys\}@physics.ox.ac.uk} \\
}

\begin{document}

\maketitle

\begin{abstract}
This paper presents a novel numerical optimisation method for infinite dimensional optimisation. The functional optimisation makes minimal assumptions about the functional and without any specific knowledge on the derivative of the functional. The algorithm has been tested on several physical systems (brachistochrone and catenary problems) and it is shown that the solutions obtained are close to the actual solutions in one thousand functional evaluations. It is also shown that for the tested cases, the new algorithm provides better convergence to the optimum value compared to the tested existing algorithms.
\end{abstract}

\section{Introduction}
Bayesian optimisation (BO) \cite{GP-machine-learning, review-ieee} is a powerful tool to optimise functions with some prior assumptions about the function. The common process in which BO is utilised is the Gaussian process (GP) \cite{GP-book}. In GP, points in the function to be optimised are assumed to be correlated with the points surrounding them. Several algorithms have been proposed using GP. Among them are GP Probability Improvement (GP-PI) \cite{GP-PI}, GP Expected Improvement (GP-EI) \cite{GP-EI}, and GP Upper Confidence Bound (GP-UCB) \cite{GP-UCB, GP-UCB-Nando}. In these algorithms, several samples of the functions are taken and then the posterior distribution of the function for the unseen positions is updated.
The algorithms then compute a number of acquisition functions for every unseen position and take a sample from the position that gives the highest acquisition function. It has been shown theoretically that GP-UCB will eventually converge to the global maximum of a function, provided that the next sampled point is the point that has the highest upper confidence bound \cite{GP-UCB}.

One major drawback of GP is that it needs an auxiliary search to look for the highest acquisition function value. Simultaneous Optimistic Optimisation (SOO) \cite{SOO} is an optimisation algorithm that does not need an auxiliary search. Several methods also take advantage of SOO and BO by combining them, such as in BaMSOO \cite{BaMSOO} and IMGPO \cite{IMGPO}. However, these algorithms only work in a moderate number of dimensions and they do not perform well in higher dimensions.

In physical applications, there are several problems that require optimisation of functionals, instead of functions. Simple examples are ``what is the trajectory of a ball to reach a point in the shortest amount of time?'' (brachistochrone), and ``what is the shape of a moving object to minimise the drag in a fluid?''. These problems can be simulated to calculate the functionals, but it is difficult to calculate the gradient of the functionals. As a result, the optimisation algorithm required is one that does not compute gradients and yet can still work for functionals. In this paper, an extension of SOO is presented to work for optimising functionals and it will be tested later to solve several physical problems.

\section{Simultaneous Optimistic Optimisation (SOO)}
Consider a function, $f: \mathcal{X} \rightarrow
\mathbb{R}$ with $\mathcal{X} \subset \mathbb{R}^D$, a bounded set in $\mathbb{R}^D$. In most cases, $\mathcal{X}$ is a hyper-rectangle with $D$ dimensions. The function has at least a global optimiser, $x^* \in \mathcal{X}$ of $f$, and for all $x \in \mathcal{X}$, it is assumed that $f(x^*) - f(x) \leq l(x^*, x) $. The function $l:\mathcal{X}\times\mathcal{X}\rightarrow \mathbb{R}$ is a semi-metric function between two points in $\mathcal{X}$.

With the properties and assumptions above, R. Munos \cite{SOO} advised optimistic optimisation algorithms to find the global optimum value of $f$. One of the algorithms, Simultaneous Optimistic Optimisation (SOO), does not even need to have the knowledge of the semi-metric function, $l$.

SOO is a tree-based search algorithm which partitions the search space, $\mathcal{X}$ into several cells. Initially, there is only one cell, which contains all the search space. The central point of the cell is evaluated. This is the root of the tree. Then the cell is divided into $K$ smaller cells and the centre points of the cells are evaluated. The smaller cells become the children of the divided cell. If a cell has a centre value larger than any cells with larger or the same size, then the cell is divided. The tree's depth is limited by a function of the number of evaluations, $h_{max}(n)$. The pseudo-code of the algorithm is shown in \ref{alg:soo}, where the $i$-th cell of depth $h$ is denoted as $(h,i)$ and the central value of the cell is $x_{h,i}$.

\begin{algorithm}
\caption{Simultaneous Optimistic Optimisation (SOO)} \label{alg:soo}
\begin{algorithmic}
\STATE $\mathcal{T} \leftarrow {(0,0)}$; $n \leftarrow 0$; and $\mathcal{L}_{\mathcal{T}} :=$ any nodes in $\mathcal{T}$ without children
\REPEAT
	\STATE $v_{max} \leftarrow 0$
	\FOR{$h \leftarrow 0$ \TO $\min(depth(\mathcal{T}), h_{max}(n))$}
		\STATE $(h,i) \leftarrow \arg\max_{(h,j)\in\mathcal{L}_{\mathcal{T}}} f(x_{h,j})$
		\IF {$v_{max} \leq f(x_{h,i})$}
			\STATE Expand $(h,i)$ to $K$ children, evaluate the values of the children, and add them to $\mathcal{T}$
			\STATE $v_{max} \leftarrow f(x_{h,i})$; $n \leftarrow n + K$
		\ENDIF
	\ENDFOR
\UNTIL any stopping conditions
\RETURN $x_n^+ \leftarrow \arg\max_{(h,i)\in\mathcal{T}} f(x_{h,i})$
\end{algorithmic}
\end{algorithm}

The regret of the algorithm after $n$ function evaluations is defined as $r_n = f(x^*) - f(x^+_n)$, where $f(x^+_n)$ is the maximum value found so far. For some functions, the regret is even exponential, i.e. $r_n = \mathcal{O}(e^{-C\sqrt{n}/D})$ for some constant, $C$. SOO is also combined with Gaussian Processes to give better empirical performance and theoretical guarantee, as seen in BaMSOO \cite{BaMSOO} and IMGPO \cite{IMGPO}.

\section{Multi-level SOO}
SOO is a powerful algorithm to optimise function with moderate dimensions (e.g. 10 dimensions). However, to optimise functionals that have infinite dimensions, some modifications are required.

Suppose that we have a 1D function, $f:\mathcal{X} \rightarrow \mathcal{F}$ where $\mathcal{X}, \mathcal{F} \subset \mathbb{R}$, which is an input of a functional, $J: \mathcal{F}^{\mathcal{X}} \rightarrow \mathbb{R}$. The problem considered here is to find a function, $f_*$, to maximise the black-box functional, $J$, i.e. $f_* = \arg\max_{f\in\mathcal{F}^{\mathcal{X}}} J[f].$

It is assumed that the end points of $f$ are fixed, i.e. at $f(x_a)$ and $f(x_b)$, where $x_a < x_b$. The functional is also assumed to be Fr\'{e}chet differentiable and $f_*$ is twice differentiable. In order to solve it numerically, the function $f(x)$ is discretised at regular spacing points and then the number of discrete points in $f(x)$ is gradually increased. The values of $f(x)$ for $x \in \mathcal{X}$ is obtained by linearly interpolating the values from the discretised points.

Let one denote $f^{(l)}(x)$ as the discretised function of $f(x)$ at $2^l-1$ points, excluding the end points, with spacing $h = (x_b - x_a) / 2^l$, and then interpolate linearly. The number of discretised points will always be counted without the end points, unless indicated. The set of all possible discretised functions is denoted as $\mathcal{F}^{\mathcal{X}(l)}$, where $\mathcal{F}^{\mathcal{X}(1)} \subset \mathcal{F}^{\mathcal{X}(2)} \subset ... \subset \mathcal{F}^{\mathcal{X}} $. The discretised optimiser function, $f_*$, is denoted by $f_*^{(l)}(x)$. The global optimiser in $\mathcal{F}^{\mathcal{X}(l)}$ is defined as
$
f_{+}^{(l)}(x) = \arg\max_{f\in\mathcal{F}^{\mathcal{X}(l)}} J[f].
$
The global optimiser in $\mathcal{F}^{\mathcal{X}(l)}$ is not necessarily the same as $f_*^{(l)}$, but $f_+^{(l)} \rightarrow f_*^{(l)}$ as $l\rightarrow\infty$.

\begin{assumption}
\label{assumption:locality-global-maximiser}
For a Fr\'{e}chet differentiable functional, $J[f]$, there exists $L > 0$ and $\eta \in (0,1)$ where $ \eta L |\!| f_* - f |\!|^2_1 \leq J[f_*] - J[f] \leq L |\!| f_* - f |\!|^2_1 $ for all $f \in \mathcal{F}^{\mathcal{X}}$.
\end{assumption}

\begin{lemma}
\label{lemma:distance-between-two-optima}
The $L_1$-distance between $f_+^{(l)}$ and $f_*$ is bounded by $|\!| f_* - f_+^{(l)} |\!|_1 \leq \alpha 4^{-l} $ for some $\alpha$.
\end{lemma}

The variable $\alpha$ depends on the intrinsic parameters of $J$ and $f_*$. From the lemma above, the search space can be decreased by $4^{-l}$ as $l$ increased, provided $\alpha$ is chosen to be large enough. Based on lemma \ref{lemma:distance-between-two-optima}, it is possible to perform optimisations using SOO as in Algorithm \ref{alg:soo} with the use of several rules in dividing cells as below, letting $p = 4$.
\begin{enumerate}
\item The initial cell is one dimensional with a width of $w$ and $l=1$. The position of the central point of the cell represents the value of a point in $x = (x_a + x_b)/2$.
\item If all dimensions in the cell have lengths less than or equal to $p^{-l}$, then new $2^l$ dimensions are added to the cell. The new dimensions represent the position of the new points in the discretised $f^{(l+1)}(x)$ which is placed in the middle positions of discretised points in $f^{(l)}(x)$. The new dimensions have widths of $p^{-l}$ each. The value of $l$ of the cell is increased by 1.
\item A cell is divided along the longest dimension into $K=3$ smaller cells as in \cite{IMGPO}. If there are more than one dimension that have the same length, then the `oldest' dimension is chosen.
\item \label{rule:dimensions-dependency} If a cell is divided along an `old' dimension, the change in position in the dimension will also change the position in the `newer' dimensions around it.
\end{enumerate}

As an example let $p = K = 3$. The first cell has $l = 1$ and only specifies $D=2^l-1=1$ points on $f$, which is $x_1^{(1)}$, the middle point between $x_a$ and $x_b$. On the first division, the cell is divided into 3 smaller cells along the first dimension. After the division, the children's dimension width is now $3^{-1}$ and it is equal to $p^{-l}$ with $p = 3$ and $l = 1$. Therefore, $2^l = 2$ new dimensions are added to the children cell which correspond to the middle points between $x_a$ and $x_1^{(1)}$, and between $x_1^{(1)}$ and $x_b$.

With $l = 2$ and $D = 2^l-1 = 3$ dimensions, let one denote the discretised points as $x_i^{(2)}$ with $i\in\{1,2,3\}$. The second dimension is the oldest dimension because it was the same dimension as the dimension with $l = 1$. The first and third dimensions are the newer dimensions. As all the dimensions now has width of $3^{-1}$, the older dimension is divided first, i.e. $i = 2$. Because an older dimension is divided, some of its children have different positions in the second dimension and thus change the positions in the first and third dimensions. If one of the children is shifted in the second dimension by $\Delta x_2^{(2)}$ with respect to its parent's position, then the first and third dimensions also shift by $\Delta x_2^{(2)}/2$, because linear interpolation is used.

For the next division, as the longest dimensions are the first and the third dimension, it can choose any dimension to be divided. As the divided dimension is the newest dimension, it does not shift the positions in other dimensions. The division process of cells is repeated until the algorithm stops.


\section{Numerical Experiments}
\label{sec:numerical-experiments}

To test the performance of the algorithm, some numerical experiments on physical cases have been performed. The algorithm is benchmarked against SOO \cite{SOO}, IMGPO \cite{IMGPO}. Those algorithms are also gradient-free optimisations, i.e. they do not need information about gradient of the functionals to optimise. As the other algorithms are designed for fixed dimensions, it is tested with 7 dimensions.

The algorithms are tested for solving brachistochrone problems. Given two points in 2D space, the algorithms need to determine the path between the points so that a bead starting from the first point can travel, in influence of gravity pointing in $-y$-direction, to the other point as fast as possible. The gravity is assumed to be 1. In the first test case, the positions of the end points are $(0,0)$ and $(1,0)$ and the beads have initial velocity $v_0=\sqrt{2/(2+\pi)}$. The minimum time required in this case is $\pi/\sqrt{2+\pi}$. In the second case, the points positions are $(0,0)$ and $(1,2/(2+\pi))$ with the initial velocity and minimum time are $v_0=2/\sqrt{2+\pi}$ and $\pi/\sqrt{4+2\pi}$, respectively. The optimum shape of the brachistochrone problems is cycloid where it can be expressed as parametric functions, $x = a(t-\sin t)$ and $y = a(1-\cos t)-a$ with some $a$. Output of the functional is the time required to travel such shape and returns infinity if it is impossible to reach the other end. Minimisation can be done by maximising the negative of the output value of the functional.

Another test case is the catenary problem, which is the minimum area when the line between two end points is rotated along the $x$-axis. The end points are located at $(-0.5,0)$ and $(0.5,0)$. The optimum shape of this problem can be expressed as a hyperbolic function, i.e. $y = a \cosh(x/a)$.

In SOO and IMGPO, a cell is divided into three smaller children along the longest dimension. The parameters that were used in IMGPO are the parameters used in \cite{IMGPO}. The search space for the algorithms is confined to $[0,1]$ for each dimension. For the new algorithm, the initial bound is set to be $[-4,4]$ and reduced by factor of 4 every time new dimensions are added.

\begin{figure}
\includegraphics[scale=0.15]{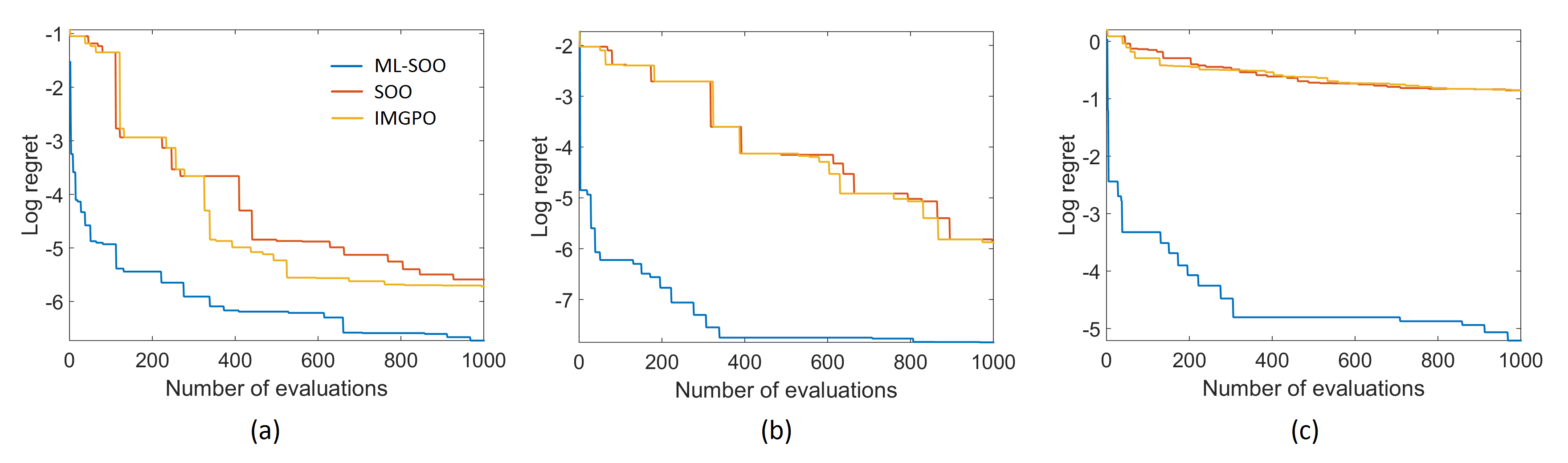}
\includegraphics[scale=0.323]{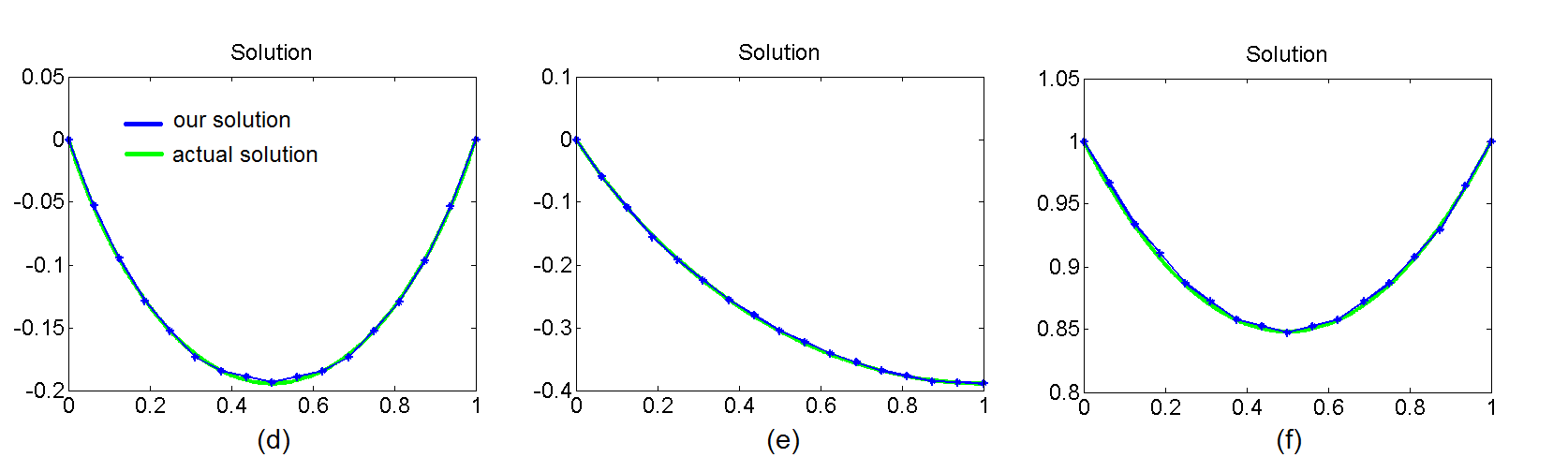}
\caption{\label{fig:results-log-regret} Top: Log regret versus number of functional evaluations using multi-level SOO, SOO, and IMGPO for (a) the first brachistochrone problem, (b) the second brachistochrone problem, and (c) the catenary problem. Bottom: Comparison of the solutions obtained by the algorithm with the actual solutions for each case. The obtained solutions have 15 discrete points, excluding the end points.}
\end{figure}

The performance of the algorithms in all test cases is shown in Figure \ref{fig:results-log-regret}(a)-(c). The solutions obtained by the new algorithm are also plotted with the optimum solution in Figure \ref{fig:results-log-regret}(d)-(f). As demonstrated in the figures, the new algorithm performs better than the other algorithms in all test cases. The solutions obtained by the new algorithm are also quite close to the optimum solutions. In the first case, the solution found by the new algorithm consists of 15 points which is already more than the dimensions tested by other algorithms. The number of discretised points in the solutions can be increased by running the algorithm for longer.

\section{Conclusion}
A novel extension of Simultaneous Optimistic Optimisation (SOO) for infinite dimensional optimisation has been presented for the first time in this paper. The objective of the algorithm is to find a curve or a 1D function to obtain the maximum value of the given functional. The algorithm has been tested against the original SOO and IMGPO with fixed dimensions to solve analytically known physical systems, such as the brachistochrone and the catenary problems. In all test cases, the new multi-level SOO gives faster convergence to the actual solution, compared to the other two algorithms. It also gives 15 discrete points in 1000 times functional evaluations.

The biggest advantage of using multi-level SOO to optimise physical systems is that it does not need the gradient of the system, just a functional to calculate the value of tested functions. It is suitable to optimise shapes in physical systems that require simulations to compute the functional values.

\subsubsection*{Acknowledgments}

The authors would like to acknowledge the supports from the plasma physics HEC Consortium EPSRC grant number EP/L000237/1, as well as the Central Laser Facility and the Computer Science Department at the Rutherford Appleton Laboratory for the use of SCARF-LEXICON computer cluster. We would also like to thank ARCHER UK National Supercomputing Service for the use of the computing service. We also wish to thank the UCLA/IST OSIRIS  consortium for the use of OSIRIS. M.F.K.  would  like  to  gratefully  thank  Indonesian  Endowment Fund for Education (LPDP) for its support. The authors gratefully acknowledge support from  OxCHEDS and P.A.N. for his William Penney Fellowship with AWE plc. M.F.K. would also like to thank Dr. Raoul Trines for the useful discussion about functionals and Ayesha Chairannisa for proof-reading the manuscript.

\end{document}